\documentclass[reqno,12pt,a4paper]{amsart}

\usepackage{mathrsfs}
\usepackage{amssymb}
\usepackage{amsfonts}
\usepackage{latexsym}
\usepackage{amsthm}

\usepackage{graphicx}
\def\lb{\label}

\newcommand{\er}[1]{\textrm{(\ref{#1})}}

\begin{document}


\renewcommand{\theequation}{\arabic{section}.\arabic{equation}}
\theoremstyle{plain}
\newtheorem{theorem}{\bf Theorem}[section]
\newtheorem{lemma}[theorem]{\bf Lemma}
\newtheorem{corollary}[theorem]{\bf Corollary}
\newtheorem{proposition}[theorem]{\bf Proposition}
\newtheorem{definition}[theorem]{\bf Definition}
\newtheorem{remark}[theorem]{\it Remark}

\def\a{\alpha}  \def\cA{{\mathcal A}}     \def\bA{{\bf A}}  \def\mA{{\mathscr A}}
\def\b{\beta}   \def\cB{{\mathcal B}}     \def\bB{{\bf B}}  \def\mB{{\mathscr B}}
\def\g{\gamma}  \def\cC{{\mathcal C}}     \def\bC{{\bf C}}  \def\mC{{\mathscr C}}
\def\G{\Gamma}  \def\cD{{\mathcal D}}     \def\bD{{\bf D}}  \def\mD{{\mathscr D}}
\def\d{\delta}  \def\cE{{\mathcal E}}     \def\bE{{\bf E}}  \def\mE{{\mathscr E}}
\def\D{\Delta}  \def\cF{{\mathcal F}}     \def\bF{{\bf F}}  \def\mF{{\mathscr F}}
\def\c{\chi}    \def\cG{{\mathcal G}}     \def\bG{{\bf G}}  \def\mG{{\mathscr G}}
\def\z{\zeta}   \def\cH{{\mathcal H}}     \def\bH{{\bf H}}  \def\mH{{\mathscr H}}
\def\e{\eta}    \def\cI{{\mathcal I}}     \def\bI{{\bf I}}  \def\mI{{\mathscr I}}
\def\p{\psi}    \def\cJ{{\mathcal J}}     \def\bJ{{\bf J}}  \def\mJ{{\mathscr J}}
\def\vT{\Theta} \def\cK{{\mathcal K}}     \def\bK{{\bf K}}  \def\mK{{\mathscr K}}
\def\k{\kappa}  \def\cL{{\mathcal L}}     \def\bL{{\bf L}}  \def\mL{{\mathscr L}}
\def\l{\lambda} \def\cM{{\mathcal M}}     \def\bM{{\bf M}}  \def\mM{{\mathscr M}}
\def\L{\Lambda} \def\cN{{\mathcal N}}     \def\bN{{\bf N}}  \def\mN{{\mathscr N}}
\def\m{\mu}     \def\cO{{\mathcal O}}     \def\bO{{\bf O}}  \def\mO{{\mathscr O}}
\def\n{\nu}     \def\cP{{\mathcal P}}     \def\bP{{\bf P}}  \def\mP{{\mathscr P}}
\def\r{\rho}    \def\cQ{{\mathcal Q}}     \def\bQ{{\bf Q}}  \def\mQ{{\mathscr Q}}
\def\s{\sigma}  \def\cR{{\mathcal R}}     \def\bR{{\bf R}}  \def\mR{{\mathscr R}}
\def\S{\Sigma}  \def\cS{{\mathcal S}}     \def\bS{{\bf S}}  \def\mS{{\mathscr S}}
\def\t{\tau}    \def\cT{{\mathcal T}}     \def\bT{{\bf T}}  \def\mT{{\mathscr T}}
\def\f{\phi}    \def\cU{{\mathcal U}}     \def\bU{{\bf U}}  \def\mU{{\mathscr U}}
\def\F{\Phi}    \def\cV{{\mathcal V}}     \def\bV{{\bf V}}  \def\mV{{\mathscr V}}
\def\P{\Psi}    \def\cW{{\mathcal W}}     \def\bW{{\bf W}}  \def\mW{{\mathscr W}}
\def\o{\omega}  \def\cX{{\mathcal X}}     \def\bX{{\bf X}}  \def\mX{{\mathscr X}}
\def\x{\xi}     \def\cY{{\mathcal Y}}     \def\bY{{\bf Y}}  \def\mY{{\mathscr Y}}
\def\X{\Xi}     \def\cZ{{\mathcal Z}}     \def\bZ{{\bf Z}}  \def\mZ{{\mathscr Z}}
\def\O{\Omega}
\def\vr{\varrho}
\def\vs{\varsigma}

\newcommand{\gA}{\mathfrak{A}}          \newcommand{\ga}{\mathfrak{a}}
\newcommand{\gB}{\mathfrak{B}}          \newcommand{\gb}{\mathfrak{b}}
\newcommand{\gC}{\mathfrak{C}}          \newcommand{\gc}{\mathfrak{c}}
\newcommand{\gD}{\mathfrak{D}}          \newcommand{\gd}{\mathfrak{d}}
\newcommand{\gE}{\mathfrak{E}}
\newcommand{\gF}{\mathfrak{F}}           \newcommand{\gf}{\mathfrak{f}}
\newcommand{\gG}{\mathfrak{G}}           
\newcommand{\gH}{\mathfrak{H}}           \newcommand{\gh}{\mathfrak{h}}
\newcommand{\gI}{\mathfrak{I}}           \newcommand{\gi}{\mathfrak{i}}
\newcommand{\gJ}{\mathfrak{J}}           \newcommand{\gj}{\mathfrak{j}}
\newcommand{\gK}{\mathfrak{K}}            \newcommand{\gk}{\mathfrak{k}}
\newcommand{\gL}{\mathfrak{L}}            \newcommand{\gl}{\mathfrak{l}}
\newcommand{\gM}{\mathfrak{M}}            \newcommand{\gm}{\mathfrak{m}}
\newcommand{\gN}{\mathfrak{N}}            \newcommand{\gn}{\mathfrak{n}}
\newcommand{\gO}{\mathfrak{O}}
\newcommand{\gP}{\mathfrak{P}}             \newcommand{\gp}{\mathfrak{p}}
\newcommand{\gQ}{\mathfrak{Q}}             \newcommand{\gq}{\mathfrak{q}}
\newcommand{\gR}{\mathfrak{R}}             \newcommand{\gr}{\mathfrak{r}}
\newcommand{\gS}{\mathfrak{S}}              \newcommand{\gs}{\mathfrak{s}}
\newcommand{\gT}{\mathfrak{T}}             \newcommand{\gt}{\mathfrak{t}}
\newcommand{\gU}{\mathfrak{U}}             \newcommand{\gu}{\mathfrak{u}}
\newcommand{\gV}{\mathfrak{V}}             \newcommand{\gv}{\mathfrak{v}}
\newcommand{\gW}{\mathfrak{W}}             \newcommand{\gw}{\mathfrak{w}}
\newcommand{\gX}{\mathfrak{X}}               \newcommand{\gx}{\mathfrak{x}}
\newcommand{\gY}{\mathfrak{Y}}              \newcommand{\gy}{\mathfrak{y}}
\newcommand{\gZ}{\mathfrak{Z}}             \newcommand{\gz}{\mathfrak{z}}

\def\ba{{\bf a}}\def\be{{\bf e}} \def\bc{{\bf c}}
\def\bm{{\bf m}}
\def\bv{{\bf v}} \def\bu{{\bf u}}

\def\be{{\bf e}} \def\bc{{\bf c}}
\def\bx{{\bf x}} \def\by{{\bf y}}
\def\bv{{\bf v}} \def\bu{{\bf u}}

\def\Om{\Omega}
\def\bbD{\pmb \Delta}
\def\mm{\mathrm m}
\def\mn{\mathrm n}

\def\ve{\varepsilon}   \def\vt{\vartheta}    \def\vp{\varphi}    \def\vk{\varkappa}

 \def\dD{{\mathbb A}}   \def\B{{\mathbb B}} \def\C{{\mathbb C}}
 \def\dD{{\mathbb D}}  \def\dE{{\mathbb E}}  \def\dG{{\mathbb G}}
 \def\dH{{\mathbb H}}
 \def\dF{{\mathbb F}}  \def\dI{{\mathbb I}}  \def\dJ{{\mathbb J}}
 \def\K{{\mathbb K}}  \def\dL{{\mathbb L}}   \def\dM{{\mathbb M}}
 \def\N{{\mathbb N}}  \def\dO{{\mathbb O}}  \def\dP{{\mathbb P}}
 \def\dQ{{\mathbb Q}}
 \def\R{{\mathbb R}}  \def\dS{{\mathbb S}}  \def\T{{\mathbb T}}
  \def\dU{{\mathbb U}}  \def\dV{{\mathbb V}} \def\dW{{\mathbb W}}
   \def\dX{{\mathbb X}} \def\dY{{\mathbb Y}} \def\Z{{\mathbb Z}}


\def\la{\leftarrow}              \def\ra{\rightarrow}      \def\Ra{\Rightarrow}
\def\ua{\uparrow}                \def\da{\downarrow}
\def\lra{\leftrightarrow}        \def\Lra{\Leftrightarrow}


\def\lt{\biggl}                  \def\rt{\biggr}
\def\ol{\overline}               \def\wt{\widetilde}
\def\no{\noindent}               \def\ti{\tilde}
\def\ul{\underline}


\let\ge\geqslant                 \let\le\leqslant
\def\lan{\langle}                \def\ran{\rangle}
\def\/{\over}                    \def\iy{\infty}
\def\sm{\setminus}               \def\es{\emptyset}
\def\ss{\subset}                 \def\ts{\times}
\def\pa{\partial}                \def\os{\oplus}
\def\om{\ominus}                 \def\ev{\equiv}
\def\iint{\int\!\!\!\int}        \def\iintt{\mathop{\int\!\!\int\!\!\dots\!\!\int}\limits}
\def\el2{\ell^{\,2}}             \def\1{1\!\!1}
\def\wh{\widehat}

\def\sh{\mathop{\mathrm{sh}}\nolimits}
\def\ch{\mathop{\mathrm{ch}}\nolimits}

\def\where{\mathop{\mathrm{where}}\nolimits}
\def\as{\mathop{\mathrm{as}}\nolimits}
\def\Area{\mathop{\mathrm{Area}}\nolimits}
\def\arg{\mathop{\mathrm{arg}}\nolimits}
\def\const{\mathop{\mathrm{const}}\nolimits}
\def\det{\mathop{\mathrm{det}}\nolimits}
\def\diag{\mathop{\mathrm{diag}}\nolimits}
\def\diam{\mathop{\mathrm{diam}}\nolimits}
\def\dim{\mathop{\mathrm{dim}}\nolimits}
\def\dist{\mathop{\mathrm{dist}}\nolimits}
\def\Im{\mathop{\mathrm{Im}}\nolimits}
\def\Iso{\mathop{\mathrm{Iso}}\nolimits}
\def\Ker{\mathop{\mathrm{Ker}}\nolimits}
\def\Lip{\mathop{\mathrm{Lip}}\nolimits}
\def\rank{\mathop{\mathrm{rank}}\limits}
\def\Ran{\mathop{\mathrm{Ran}}\nolimits}
\def\Re{\mathop{\mathrm{Re}}\nolimits}
\def\Res{\mathop{\mathrm{Res}}\nolimits}
\def\res{\mathop{\mathrm{res}}\limits}
\def\sign{\mathop{\mathrm{sign}}\nolimits}
\def\span{\mathop{\mathrm{span}}\nolimits}
\def\supp{\mathop{\mathrm{supp}}\nolimits}
\def\Tr{\mathop{\mathrm{Tr}}\nolimits}
\def\BBox{\hspace{1mm}\vrule height6pt width5.5pt depth0pt \hspace{6pt}}


\newcommand\nh[2]{\widehat{#1}\vphantom{#1}^{(#2)}}
\def\dia{\diamond}

\def\Oplus{\bigoplus\nolimits}



\def\qqq{\qquad}
\def\qq{\quad}
\let\ge\geqslant
\let\le\leqslant
\let\geq\geqslant
\let\leq\leqslant
\newcommand{\ca}{\begin{cases}}
\newcommand{\ac}{\end{cases}}
\newcommand{\ma}{\begin{pmatrix}}
\newcommand{\am}{\end{pmatrix}}
\renewcommand{\[}{\begin{equation}}
\renewcommand{\]}{\end{equation}}
\def\bu{\bullet}


\title[{The number of eigenvalues of discrete  Hamiltonian periodic in time}]
 {The number of eigenvalues of  discrete  Hamiltonian periodic in time}

\author[Evgeny, L. Korotyaev]{Evgeny, L. Korotyaev}
\address{E. Korotyaev, Depart. of Math. Analysis, Saint-Petersburg State University,
Universitetskaya nab. 7/9, St. Petersburg, 199034, Russia, \
korotyaev@gmail.com, \ e.korotyaev@spbu.ru}

\date{\today}

\begin{abstract}
\no  We consider time periodic Hamiltonian on periodic graphs and
estimate the number of its quasi-energy eigenvalues  on the finite
interval.

\end{abstract}


\maketitle


\section {Introduction}

We consider  Schr\"odinger equations on a periodic graph $\cG$ given
by
\[
\lb{eg1xx}
\begin{aligned}
{\tfrac{ d }{d t}}u(t)= -i h(t) u(t),\qq h(t)=\D+\gp+ V(t),
\end{aligned}
\]
where $h(t) $ is the Hamiltonian,  $\t$-periodic in time $t$; $\D$
is the discrete   Laplacian on $\G$ and $\gp:\cV\to \R $ is a
bounded electric potential. We assume that the potential $V(t)$ on
$\cG$ is $\t$ -periodic in time and is decaying on $\cG$.

There are a lot of applications of different periodic media, e.g.,
nanomedia, in physics, chemistry and engineering, see, e.g.,
\cite{NG04}, two authors of this survey, Novoselov and Geim, won the
Nobel Prize for discovering graphene. In order to study properties
of such media one uses their approximations by periodic graphs. In
these cases we need to investigate various properties of
Schr\"odinger operators on graphs, and, in particular, periodic in
time.  One of the interesting problems is  Schr\"odinger operators
with  decaying potentials, periodic in time, on periodic graphs. As
far as the author knows, there are no results about eigenvalues for
Schr\"odinger operators with  decaying potentials, periodic in time,
on periodic graphs.

\subsection{Schr\"odinger operators on periodic graphs}

Let $\cG=(\cV,\cE)$ be a connected infinite graph, possibly  having
loops and multiple edges and embedded into the space $\R^d$. Here
$\cV$  is the set of its vertices and $\cE$ is the set of its
unoriented edges.
An edge starting at a vertex $x$ and ending at a vertex $y$ from
$\cV$ will be denoted as the ordered pair $(x,y)\in\cE$ and is said
to be \emph{incident} to the vertices.
Vertices $x,y\in\cV$ will be called \emph{adjacent} and denoted by
$x\sim y$, if $(x,y)\in \cE$. We define the degree $\vk_x$ of the
vertex $x\in\cV$ as the number of all edges from $\cE$, starting at
$x$. Let $\G$ be a lattice of rank $d$ in $\R^d$ with a basis
$\{a_1,\ldots,a_d\}$, i.e.,
$$
\G=\Big\{a : a=\sum_{s=1}^dn_sa_s, \; (n_s)_{s=1}^d\in\Z^d\Big\}.
$$
We introduce the equivalence relation on $\R^d$:
$$
\bx\equiv \by \; (\hspace{-4mm}\mod \G) \qq\Leftrightarrow\qq
\bx-\by\in\G \qqq \forall\, \bx,\by\in\R^d.
$$

We consider \emph{locally finite $\G$-periodic graphs} $\cG$, i.e.,
graphs satisfying the following conditions:

  $\bu $ $\cG=\cG+a$ for any $a\in\G$ and the quotient graph  $\cG_*=\cG/\G$ is finite.

The basis $a_1,\ldots,a_d$ of the lattice $\G$ is called the {\it
periods}  of $\cG$. We also call the quotient graph $\cG_*=\cG/\G$
the \emph{fundamental graph} of the periodic graph $\cG$. The
fundamental graph $\cG_*$ is a graph on the $d$-dimensional torus
$\R^d/\G$. The graph $\cG_*=(\cV_*,\cE_*)$ has the vertex set
$\cV_*=\cV/\G$.
Introduce the space $\ell^p(\cV), p\ge 1$  of sequences
$f=(f_x)_{x\in \cV}$ equipped with the norm given by
$$
\begin{aligned}
\|f\|_{p}=\|f\|_{\ell^p(\cV)} =\big(\sum_{x\in
\cV}|f_x|^p\big)^{1\/p},\qq  \ p\in [1,\iy),
\end{aligned}
$$
and let $\|f\|_{\iy}=\|f\|_{\ell^\iy(\cV)} =\sup_{x\in \cV}|f_x|$
and $\|f\|=\|f\|_{2}$ for shortness.

We consider the Schr\"odinger operator $h_o$ acting on the Hilbert
space $\ell^2(\cV)$ and given by
\[
\lb{Sh} h_o=\D+\gp,
\]
where the potential $\gp$ is real valued, bounded and satisfies for
all $(x,a)\in\cV\ts\G$:
\[
\lb{Pot} (\gp f)_x=\gp_xf_x, \qqq \gp_{x+a}=\gp_x,
\]
and $\D$ is the \emph{combinatorial Laplacian} on $f\in\ell^2(\cV)$
defined by
\[
\lb{DLO} (\D f)_x=\sum_{ y\sim x}\big(f_x-f_y\big), \qqq x\in\cV.
\]
The sum in \er{DLO} is taken over all edges starting at the vertex
$x$. It is known   that $\D$ is self-adjoint, the point 0 belongs to
its spectrum $\s(\D)$ containing in $[0,s_+]$ for some $s_+\le
2\vk_+$, i.e.,
\[
\lb{bf} 0\in\s(\D)\ss [0,s_+]\subset[0,2\vk_+],\qqq
\textrm{where}\qqq \vk_+=\max_{x\in\cV}\vk_x<\iy,
\]
see, e.g., \cite{MW89} and $\vk_x$ is the degree of the vertex $x$,
and $s_+$ is the top of the spectrum $\s(h_o)$.

We briefly describe the spectrum of the Schr\"odin\-ger operator
$h_o$ (for more details see, e.g., \cite{HS99} or \cite{KS14}). Let
$\#M$ denote the number of elements in a set $M$. The operator $h_o$
on $\ell^2(V)$ has the standard decomposition into a constant fiber
direct integral by
\[
\lb{raz}
\ell^2(V)={1\/(2\pi)^d}\int^\oplus_{\T^d}\ell^2(V_\ast)\,d\vt ,\qqq
Uh_o U^{-1}={1\/(2\pi)^d}\int^\oplus_{\T^d}h_o(\vt)d\vt,
\]
$\T^d=\R^d/(2\pi\Z)^d$, for some unitary operator $U$. Here
$\ell^2(V_\ast)=\C^\nu$ is the fiber space, $\nu=\#V_*$ and
$h_o(\vt)$ is the fiber operator acting on $\ell^2(V_\ast)$. Each
fiber operator $h_o(k)$, $k\in\T^{d}$, acts on the space
$\ell^2(\cV_*)=\C^\n$, $\n=\#\cV_*$, and has $\n$ real eigenvalues
$\gb_j(\vt)$, $j=1,\ldots,\n$, labeled in non-decreasing order by
\[
\label{eq.3H}
\gb_{1}(\vt)\leq\gb_{2}(\vt)\leq\ldots\leq\gb_{\nu}(\vt), \qqq
\forall\,\vt\in\T^{d},
\]
counting multiplicities.  Each $\gb_j(\cdot)$ is  a real and
piecewise analytic function  on the torus $\T^{d}$ and creates the
\emph{spectral band} $\s_j(H)$ by
\[\lb{ban.1H}
\s_j(h_o)=[\gb_j^-,\gb_j^+]=\gb_j(\T^{d}),\qqq j\in\N_\n, \qqq
\N_\n=\{1,\ldots,\n\}.
\]
Some of $\gb_j(\cdot)$ may be constant, i.e.,
$\gb_j(\cdot)=\gb_j^o=\const$,  on some subset $\cB$ of $\T^d$ of
positive Lebesgue measure. In this case the Schr\"odinger operator
$h_o$ on $\cG$ has the eigenvalue $\gb_j^o$ of infinite
multiplicity. We call $\{\gb_j^o\}$ a \emph{flat band}. Thus, the
spectrum of the Schr\"odinger operator $h_o$ on the periodic graph
$\cG$ has the form
\[
\lb{specH}
 \s(h_o)=\bigcup_{k\in\T^d}\s\big(h_o(k)\big)=
\bigcup_{j=1}^{\nu}\s_j(h_o)=\s_{ac}(h_o)\cup \s_{fb}(h_o),
\]
where $\s_{ac}(h_o)$ is the absolutely continuous spectrum, which is
a union of non-degenerate bands from \er{ban.1H}, and $\s_{fb}(h_o)$
is the set of all flat bands (eigenvalues of infinite multiplicity).

{\bf Definition.} {\it The top $s_+$ of the spectrum $\s(h_o)\ss
[0,s_+]$  (i.e., the right point of the spectrum) is said to be
regular if the set $\{\vt\in \T^d: \gb_\n(\vt)=s_+\}$ is finite
number of points $\vt_j, j=1,...,n_R$ and each point satisfies
$$
\gb_\n(\vt)=s_+-(\gB_j a,a)+O(|a|^3), \qqq \as \ a=\vt-\vt_j\to
0,\qq j=1,..,n_R,
$$
for some positive matrix $\gB_j>0$.}

The regularity of the left boundary of the spectrum is defined
similarly. Korotyaev and Saburova (see Theorem 1.2 in \cite{KS16}),
show that the lower bound of the spectrum $\s(h_o)$ of the operator
$h_o$ is regular.

 We recall results of Korotyaev and Sloushch \cite{KS20a} about
estimates for the discrete spectrum of Schr\"odinger operators on
discrete periodic graphs. In order to formulate these results   we
need some definitions. A graph is called \emph{bipartite} if its
vertex set is divided into two disjoint sets (called \emph{parts} of
the graph) such that each edge connects vertices from distinct sets
(see p.105 in \cite{Or62}). Examples of bipartite graphs are the
cubic lattice and the hexagonal lattice. The face-centered cubic
lattice  is non-bipartite. Due to \cite{KS16}   the top $s_+$ of the
spectrum $\s(h_o)\ss [0,s_+]$ (i.e., the right point of the
spectrum) is regular for a bipartite graph. In general, it is not
true.

We have the decomposition for $z\in \R$:
$$
z=(z)_+-(z)_-, \qqq (z)_\pm =(|z|\pm z)/2.
$$

\begin{theorem}
\lb{Tks1}

Let $ h= h_o+\gq$ on a periodic graph, where the operator
$h_o=\D+\gp$ and the potential $\gp$  is $\G-$ periodic. Assume that
$\gq\in \ell^{d\/2}(\cV), d\ge 3$. Then the number $\gn_-(\gq,\m)$
of eigenvalues of $h$ (counted according to their multiplicity)
below $\m<0$ is finite and satisfies
\[
\lb{n-} \gn_-(\gq,\m)\le C_{ks}^- \sum
_{x\in\cV}(\gq_x-\m)_-^{d\/2}, \qq \m<0.
\]
Let in addition the top $s_+$ of the spectrum $\s(h_o)\ss [0,s_+]$
is regular (or the graph is bipartite). Then the number
$\gn_+(\gq,\m)$ of eigenvalues of $h$ (counted according to their
multiplicity) above $\m>s_+$ is finite and satisfies
\[
\lb{n+} \gn_+(\gq,\m)\le C_{ks}^+ \sum
_{x\in\cV}(\gq_x-\m)_+^{d\/2}, \qq \m>s_+.
\]
Here the constants $C_{ks}^\pm$ depend on the graph $\cG$ and the
potential $\gp$ only.

\end{theorem}

\no {\bf Remark.} In the case $d=1$ in Section 4 the discrete analogue
of Bargmann's bound from \cite{HS02} will be used. I do not know
similar estimates for $d=2$.

Below we need the  Rellich type Theorem \ref{TKo}. It is
generalization of results of  Isozaki and  Maioka \cite{IM14} and
Vesalainen \cite{V14}.

Define a half space $\dY_j(z)$ and a rectangular cuboid (or simply a
cuboid) $\dQ(\gm)$ by
\[
\begin{aligned}
\lb{XP} \dY_j(z)=\{x\in \Z^d:  x_j\ge z\},\qqq z\in \N,\ \ j\in
\N_d,
\\
\dQ(\gm)=\{x\in \Z^d:  0\le x_j\le \gm_j\},\qqq \gm=(\gm_j)\in \N^d.
\end{aligned}
\]

\begin{theorem}
\lb{TKo}

 Let $\l\in  (0, 2d)$ and $u: Z^d\to \C $  satisfy
$\sum _{|x|<r}|u_x|^2=o(r)$ \ as \ $ r\to \iy$. Let $ (\D-\l)u=f$,
where $e^{t|x|}f_x\in \ell^2(\Z^d)$ for all $t>0$.

i) If $f\big|_{\dY_d(z)}=0$ for some $z\in \Z$, then we have
$u\big|_{\dY_d(z)}=0$.

ii) If $E=\Z^d\sm \dQ(\gm)$ for some $\gm\in \N^d$ and
$f\big|_{E}=0$, then we have $u\big|_{E}=0$.

\end{theorem}

\subsection { Main results} Consider the Schr\"odinger equation   $
{\tfrac{ d }{d t}}y(t)= -i h(t) y(t)$ on the graph $\cV$, where
$h(t)=h_o+ V(t) $ is the Hamiltonian, $\t$-periodic in time $t$ and
$h_o=\D+\gp$ is the operator  defined by \er{Sh}.
 Here the potential $V_x(t), (t,x)\in \R\ts \cV$ is real
 $\t$-periodic in time
  and satisfies
\[
\lb{Vu} \sup_{t\in \T_\t}|V(t)|\le u,\qqq   \where \qq u\in
\ell^{d}(\cV).
\]
In Section 1, we assume that $d\ge 3$. The case $d=1$ is discussed
in Section 4.

Let $\pa=-i\frac{\pa}{\pa t}$ be the  self-adjoint operator in
$L^2(\T_\t)$. We also denote $\pa=-i\frac{\pa}{\pa t}$ the
corresponding operator in $\wt \cH=L^2(\T, \ell^2(\Z^d))$ with the
natural domain $\mD=\mD(\pa)$. We use the notation  $\lan A(t)\ran$
to indicate multiplication  by $A(t)$ on the space
$\cH=\ell^{2}(\Z^d)$. Introduce the {\it quasienergy} operators $\wt
h_o$ and $\wt h$ on $\wt\cH$ by
$$
\wt h_o=\pa+h_o,\qqq  \wt h=\wt h_o+V,\qq V=\lan V_x(t) \ran.
$$
 The spectrum of
$\wt h_o$ has the form
\[
\lb{sp1} \s(\wt h_o)=\s_{ac}(\wt h_o)=\bigcup_{n\in \Z} \s_n,\qqq
\s_n= \s(h_o+\o n).
\]
Note that if $\o >s_+$, then the spectrum of $\wt h_o$ has the band
structure with the sets $ \s_n$ separated by gaps.
  The spectrum of $\wt h$ has the form
$$
\s(\wt h)=\s_{ac}(\wt h)\cup \s_{pp}(\wt h),\qqq \s_{ac}(\wt
h)=\s(\wt h_o).
$$
 Introduce the  resolvent  $R(\l)=(\wt h-\l)^{-1}, \l\in \C_\pm$.
We have the  identity
\[
\lb{i1} \lan e^{it\o}\ran R(\l)\lan e^{-it\o}\ran  =R(\l+\o),\qq
\forall \l\in \C_\pm.
\]
It means that the spectrum of $\wt h$ (and $\wt h_o$)  is
$\o$--periodic. Thus it is sufficient to study eigenvalues of $\wt
h$ on the set $(-\o,0]$. The operator $\wt h$ has $N\le \iy$
eigenvalues $\{\l_j, j=1,....,N\}$ on the interval   $(-\o,0]$. Here
and below each eigenvalue is counted according to its  multiplicity.

\smallskip


\smallskip

\bigskip

$\bu $ Let $\gn(q,s)$ be the number of eigenvalues of $h_o-q$ on the
set $(-\iy,s), s\le0$ counted according to their multiplicity, where
$q$ is a potential.

$\bu $ Let $\wt \gn(A)$ be the number of eigenvalues of $\wt h$ on
the segment $A\ss \R$, counted according to its multiplicity.

Firstly, we consider the case $\o>s_+$, when there is the main gap
$\g\ss (-\o,0]$ defined by
\[
\lb{h}  \g=(-s_\g,0)\ne  \es, \qq \where \qq s_\g=\o-s_+>0.
\]
We estimate the number $\wt \gn(\g)$of eigenvalues of $\wt h$ on the
main gap $\g$.

\begin{theorem}
\lb{T1}

 Let $\wt h=\wt h_o+V$ on
$L^2(\T_\t,\ell^2(\cV)), d\ge 3$ and let $\o>s_+$.

\no i) Let potentials $V$ and $q_-=\gc u^2+v_-$, where
$\gc={16\o\/s_e^2}$ satisfy
\[
\lb{V1}  \sup_{t\in \T_\t}|V(t)|\le u,\qq  \sup_{t\in \T_\t}
V(t)_-\le v_-,\qqq \where \qq u^2, v_-\in \ell^p(\cV),\ \  p={d\/2}.
\]
 Then the operator $\wt h$ has $\wt \gn(\g_-)$ eigenvalues in
$\g_-=[-{s_\g\/2},0)$ such that
\[
\lb{11}
\begin{aligned}
 \wt \gn(\g_-)\le
 {9\o\/s_\g} C_{ks}^- \rt(1+2{\|q_-\|_\iy\/\o}\rt)\|q_-\|_p^p.
\end{aligned}
\]
ii) Let the top $s_+$ of the spectrum $\s(h_o)\ss [0,s_+]$ is
regular (or the graph is bipartite) and let
\[
\lb{V2} sup_{t\in \T_\t}|V(t)|\le u,\qq
  \sup_{t\in \T_\t}V(t)_+\le
v_+,\qqq \where \qq u^2, v_+\in \ell^{d\/2}(\cV).
\]
Let a  potential $q_+=\gc u^2+v_+$ and $\g_+=(-s_\g,-{s_\g\/2}]$.
Then
\[
\lb{12}
\begin{aligned}
\wt \gn(\g_+)\le {9\o\/s_\g} C_{ks}^+
\rt(1+2{\|q_+\|_\iy\/\o}\rt)\|q_+\|_p^p.
\end{aligned}
\]
iii) Let $V$ satisfy \er{V1} and \er{V2}. Then the number of
eigenvalues $\wt \gn(\g)$ of the operator $\wt h$ in the main gap
$\g=(-s_\g,0)$ is finite and satisfies
\[
\lb{12xx}
\begin{aligned}
 \wt \gn(\g)\le  C^+\|q_+\|_p^p+C^-\|q_-\|_p^p,
\end{aligned}
\]
where $C^\pm={9\o\/s_\g} C_{ks}^\pm \rt(1+2{\|q_\pm\|_\iy\/\o}\rt)$.

\end{theorem}

We discuss estimates for $\o$ large  enough.

\begin{theorem}
\lb{T2}

 Let $\wt h=\wt h_o+V$ on
$L^2(\T_\t,\ell^2(\cV))$. Let $\d_{\pm}=\sup_{(t,x)\in \T_\t\ts \cV}
V_x(t)_\pm<{s_\g\/2}$.

\no i) Let $V$ satisfy \er{V1} and let a  potential $q_-=c
u^2+v_-\ge 0, c>{4\o\/s_\g^2}$ and $\d_+c\le \sqrt3$ and $\o$
satisfy
\[
\lb{14} \o\ge s_++\d_-+\sqrt{\d_-(2s_++\d_-)}.
\]
 Then the operator $\wt h$ has $\wt \gn(\g_-)$ eigenvalues on
$\g_-=[\d_-,0)$ such that
\[
\lb{15}
\begin{aligned}
 \wt \gn(\g_-)\le  \gn_-(q,0)\le C_{ks}^-\|q_{_-}\|_p^p.
\end{aligned}
\]
ii) Let the top $s_+$ of the spectrum $\s(h_o)\ss [0,s_+]$ is
regular (or the graph is bipartite). Let $V$ satisfy \er{V2} and let
a potential $q_+=c u^2+v_+, c>{4\o\/s_\g^2}$ and $\d_-c\le \sqrt3$
and $\o$ satisfy
\[
\lb{16} \o\ge s_++\d_++\sqrt{\d_+(2s_++\d_+)}.
\]
 Then the operator $\wt h$ has $\wt \gn(\g_+)$ eigenvalues on
$\g_+=(-s_\g,\d_--s_\g]$, such that
\[
\lb{17}
\begin{aligned}
 \wt \gn(\g_+)\le  \gn_+(q,s_+)\le C_{ks}^+\|q_{_+}\|_p^p.
\end{aligned}
\]
 iii) Let conditions of i)-ii) hold true. Then the
operator $\wt h$ has $\wt \gn(\g)$ eigenvalues on the main gap
$\g=(-s_\g,0)$ such that
\[
\lb{18}
\begin{aligned}
 \wt \gn(\g)\le C_{ks}^+\|q_{_+}\|_p^p+C_{ks}^-\|q_{_-}\|_p^p.
\end{aligned}
\]

\end{theorem}

\subsection{The imbedding eigenvalues}
We consider the imbedding eigenvalues for unperturbed operators
$h_o=\D$ on the lattice $\Z^d, d\ge 1$, where  $\s(\D)=[0,s_+],
s_+=4d$. There are two cases:

 If $\o>s_+$, then there exist the main gap
$\g=(s_b,0)$ and the main band $\s_b=(-\o, -s_b)$.

If $\o\le s_+$, then we have $\s(\wt h_o)=\R$.

\no The main spectral interval has the form $(-\o,0)=\s_e\cup
\s_b\cup \{s_e\}$, where
\[
\gr=[s_+/\o]\in \Z_+,\qqq  s_e=\gr \o-s_+,\qq \s_b=(-\o, -s_e),\qq
\s_e=(-s_e,0).
\]
\no If $\gr=0$, then there are gaps in the spectrum of $\wt h_o$ and
$\s_e=(-s_e,0)=\g$ is the main  gap.

\no If $\gr\in \N$, then there are no gaps in the spectrum of $\wt
h_o$ and $\s_b=(-\o,0), \s_e=\es$.

We discuss the number of eigenvalues imbedding on the main spectral
interval $(-\o,0)$.

\begin{theorem}
\lb{T3}

Let $\wt h=\wt h_o+V$ on $L^2(\T_\t,\ell^2(\Z^d))$. Assume that each
$V(t), t\in \T_\t$ is finitely supported, satisfies \er{V1}, \er{V2}
and for some cuboid $\dQ_d(\gm), \ \gm\in \N^d$:
\[
\lb{Vc} \supp V_x(t)\ss \dQ(\gm) \ss \Z^d,\qq \forall \ t\in \T_\t,
\]
 Let
 $q_{\a,\pm}=\gc_\a u^2+v_\pm, \a\in \{b,e\}$, where
$\gc_e={16\o(\gr+1)\/s_e^2}$ and  $\gc_b={4\o(\gr+2)\/s_e^2}$ and
$\gr=[s_+/\o]\ge 0$.

i) Let $\gr\ge  0$. Then the operator $\wt h$ has $\wt \gn(\s_b)$
eigenvalues on $\s_b=(-\o, -s_e)$ such that
\[
\lb{19} \wt \gn(\s_b)\le  C_b \rt(\gr+2+
2{\|q_{b,\pm}\|_\iy\/\o}\rt)C_{ks}^\pm\|q_{b,\pm}\|_p^p,
\]
where $C_b=1+  {4\o(\o-s_e)(\gr+2)\/s_e^2}$.

ii) Let  $\gr\ge 1$. Then the operator $\wt h$ has $\wt \gn(\s_e)$
eigenvalues on $\s_e=(-s_e, 0)$ such that
\[
\lb{Nne} \wt \gn(\s_e)\le C^+\|q_{e,+}\|_p^p+C^-\|q_{e,-}\|_p^p,
\]
where $ C^\pm=C_e C_{ks}^\pm\rt(\gr+1+ 2{\|q_{e,\pm}\|_\iy\/\o}\rt)$
and $C_e=1+{8\o(\gr+1)\/s_e}$.
\end{theorem}

\no {\bf Remark.} 1) We consider the case $\Z^d, d\ge 3$, the case $d=1$
is discussed in Section 4.

\subsection{Historical review}
For time-periodic Hamiltonians many papers have been devoted to
scattering  mainly for self-adjoint operators $h(t)=h_o + V(t,x),
h_o=-\D$ on $\R^d, d\ge 1$, and to the spectral analysis of the
corresponding monodromy operator. Zel'dovich \cite{Z73} and Howland
\cite{H79} reduced
 the problem with a time-dependent Hamiltonian $h(t)$ to a problem with a
time-independent Hamiltonian $\wt h=-i{\pa\/\pa t}+h(t)$ by
introducing an additional time coordinate. Completeness of the wave
operators for $\wt h, \wt h_o=-i{\pa\/\pa t}+h_o$ was established by
Yajima \cite{Y77}. It was shown by Howland \cite{H79}, that $\wt h$
has no singular continuous spectrum. Korotyaev \cite{K84} proved in
that the total number of embedded eigenvalues on the interval
$[0,\o]$, counting multiplicity, is finite. The case of Schrodinger
operators with time-periodic electric and homogeneous magnetic field
was discussed in \cite{K80}, \cite{K89}, \cite{Y82}, see also recent
papers \cite{AK19}, \cite{AKS10}, \cite{AK16}, \cite{Ka19},
\cite{M00}. Moreover, scattering for three body systems was
considered in \cite{K85}, see also \cite{MS04}.

Now we discuss stationary case of Schr\"odinger operators on the
cubic lattice $\Z^d, d\ge 2$, when potentials are real and  do not
depend on time.  For Schr\"odinger operators with decaying
potentials on the lattice $\Z^d$, Boutet de Monvel and Sahbani
\cite{BS99} used Mourre's method to prove completeness of the wave
operators, absence of singular continuous spectrum and local
finiteness of eigenvalues away from threshold energies. Isozaki and
Korotyaev \cite{IK12} studied the direct and the inverse scattering
problem as well as trace formulas. Korotyaev and Moller \cite{KM19}
discussed the spectral theory for potentials $V\in \ell^p, p>1$.
Isozaki and  Morioka \cite{IM14} and Vesalainen \cite{V14} proved
that the point-spectrum of $H$ on the interval $(0,2d)$ is absent,
see also \cite{AIM16}. An upper bound on the number of discrete
eigenvalues in terms of some norm of potentials was given by
Korotyaev and Sloushch \cite{KS20}, Rozenblum and
Solomyak\cite{RoS09}. For closely related problems, we mention that
Parra and Richard \cite{PR18} reproved the results from \cite{BS99}
for periodic graphs.  Finally, scattering on periodic metric graphs
 was considered by Korotyaev and Saburova
\cite{KS20}.


\section {Eigenvalues in the main gap for $\o>s_+$}
\setcounter{equation}{0}

\subsection {Preliminaries}

 Define operators
\[
\lb{ft} \pa=\pa_+-\pa_-,\qq \pa_\pm=|\pa|\1_{\pm \pa},\qqq \where
\qq \1_s=\ca 1, s\ge0 \\ 0, s<0 \ac.
\]
Let the main spectral interval $(-\o,0)=\s_e\cup \s_b\cup \{s_e\}$,
where
\[
\s_b=(-\o, -s_e),\qq \s_e=(-s_e,0) \qqq  s_e=(\gr+1) \o-s_+,\qq
\gr=[s_+/\o]\ge 0,
\]
where $[z]$ is the integer part of $z\ge 0$. Define a projector
$P(\gr)$ in $\ell^2(\Z)$ by
\[
\lb{dP} (P(\gr)f)_n=\d_{n,j}f_n, \qq j\in \N_r^-=\{-1,..,-r\}, \qq
\forall \ n\in \Z.
\]

\begin{lemma}
\lb{Tr1} Let $\gr=[s_+/\o]\ge 0$ and  $\gc_b={4\o(\gr+2)\/s_e^2}$
and $\gc_e={16\o(\gr+1)\/s_e^2}$. Then
\[
\lb{r1} \|\pa_-^{1\/2}P^\bot(\gr) R_o(\m)\|^2\le c_e/4, \qq if \
\m\in \s_e^-=[-s_e/2,0),
\]
\[
\lb{r2} \|\pa_-^{1\/2}P^\bot(\gr+1) R_o(\m)\|^2\le c_b/4, \qq if \
\m\in \s_b=(-\o,-s_e).
\]
\end{lemma}
\no{\bf Proof.} We have identities
$$
\pa_-^{1\/2} R_o(\m)P^\bot(\gr) f= \rt( {|\o n|^{1\/2}\/\o
n+h_o-\m}f_n\rt)_{n\le -r-1},\qq
  \where \qq f=(f_n)_{n\in \Z}\in \ell^2(\cH).
$$
and the estimate
$$
\|{|\o n|^{1\/2}\/\o n+h_o-\m}\|\le \z_n,\qqq \z_n={|\o
n|^{1\/2}\/\o |n|-s_+-|\m|}, \qq n\le -\gr-1.
$$
 For $n\le -\gr-1$ and $\m\in \s_e^-$ we have

\no $\bu $ if $n= -\gr-1$,  then
$$
\z_{-\gr-1}={(\gr+1)\o)^{1\/2}\/(\gr+1)\o-s_+-|\m|}=
{((\gr+1)\o)^{1\/2}\/s_e-|\m|}\le {((\gr+1)\o)^{1\/2}\/s_e/2}=
{\sqrt {c_e}\/2},
$$
$\bu $ if $n\le -\gr-2$, then
$$
\z_n={\sqrt {\o n}\/\o n-s_+-|\m|}<\z_{-r-1}<{\sqrt {c_e}\/2}.
$$
Collecting these estimates we obtain   \er{r1}. The proof of \er{r2}
is similar. \BBox

Let $\p_j(x,t), j=1,..,N$ be orthonormal eigenfunctions of $\wt
h=h_o+V$ such that
\[
\lb{ei1} \wt h\p_j=\l_j \p_j,\qq \l_j\in (-\o,0).
\]
Let $\gM(A)$ be linear span of these  eigenvectors $\p_j$  for
$\l_j\in A$ for some segment  $A\ss (-\o,0)$.

\begin{theorem}
\lb{Tei1}

 Let $\wt h=\wt h_o+\lan V(t) \ran$, where $V$ satisfies \er{V1} and let
 $\o>s_+$, i.e., $\gr=0$.

\no i) Let $\m\in \g_-=[-{s_\g\/2},0)$ and $I_\d=[\m-\d,\m)\ss \g_-$
for some $\d\le {1\/c_e}$, where $c_e={16\o(\gr+1)\/s_e^2}$. Let a
potential $q_-=c_e u^2+v_-$.
Then the operator $\wt h$ has $\wt \gn(I_\d)$ eigenvalues, counted
with its multiplicity, on $I_\d$ such that
\[
\lb{N1} \wt \gn(I_\d)\le \sum _{n\in \Z}\gn_-(q_-,\m-\o |n|).
\]
ii) Let $I_\d=[-\d,0)$, where $\d={1\/c}$ for some
$c>{4\o\/s_\g^2}$. Let a potential $q_o=c u^2+v_-\ge 0$. Then
\[
\lb{N2} \wt \gn(I_\d)\le \sum _{n\in\Z}\gn(q_o,-\o |n|).
\]
\end{theorem}

\no {\bf Proof.} i) Let $\p=\sum g_j\p_j$, where $\p, \p_{j}\in
\gM(I_\d)$ and $g_j\in \C$. Define the operators $
T_o=\pa_-+\pa_++h_o$. We have the identities
\[
\lb{11e}
\begin{aligned}
((T_o-\m)\p,\p)=((\pa_-+\pa_++h_o-\m)\p,\p)
 =  2(\pa_-\p,\p)-(V\p,\p)+((\wt h-\m)\p,\p).
\end{aligned}
\]
Here $((\wt h-\m)\p,\p)\le 0$ and we need to estimate
$(\pa_-\p,\p)=\|\pa_-^{1\/2}\p\|^2$ in terms of $\|V\p\|$ and $((\wt
h-\m)^2\p,\p)$. Then using the identity
$$
\begin{aligned}
\p =R_o(\m)(\wt h_o-\m)\p=R_o(\m)\big((\wt h-\m)-V\big)\p,
\end{aligned}
$$
and the estimate \er{r1} at $\gr=0$,   we obtain
$$
\begin{aligned}
2\|\pa_-^{1\/2}\p\|^2= 2\|\pa_-^{1\/2} R_o(\m)\big((\wt
h-\m)-V\big)\p\|^2\le {c_e\/2} \|\big((\wt h-\m)-V\big)\p\|^2
\\
\le c_e\big(\|(\wt h-\m)\p\|^2+\|V\p\|^2 \big).
\end{aligned}
$$
Thus jointly with \er{11e} and  $c_e V^2-V\le q=c_e u^2+v_-$ we
obtain
\[
\lb{qw1z}
\begin{aligned}
((T_o-\m)\p,\p)=2(\pa_-\p,\p)-(V\p,\p)+((\wt h-\m)\p,\p)
\\
\le c_e\Big(\|V\p\|^2+\|(\wt h-\m)\p\|^2\Big)-(V\p,\p)+((\wt
h-\m)\p,\p)
\\
\le ((c_e V^2-V)\p,\p)\le (q_-\p,\p),
\end{aligned}
\]
where we have used $c_e\|(\wt h-\m)\p\|^2+((\wt h-\m)\p,\p)\le 0$,
since $c_e(\l_j-\m)^2+(\l_j-\m)\le 0$.  We rewrite \er{qw1z} in the
form
\[
\lb{23x} ((T_o-\m)\p,\p)\le (q_-\p,\p).
\]
where $(\F T_o \F^* f)_n=(\o |n|+h_o-q)f_n$ and $f_n=(\F f)_n $ for
$ n\in \Z$. From this and \er{n-}, \er{23x}  we obtain \er{N1}.

ii) Let $-\d<\ve<0$ and $\ve\to 0$. In \er{N1} we have estimated
$\wt \gn(I(\ve))$, where $I(\ve)=[-\d, \ve)$ in terms of
$\gn(q(\ve),\ve-\o n)$ for a finite $n\le 0$, where the potential
$q(\ve)=C_\ve u^2+v\ge 0$ and $C_\ve={4\o\/(s_\g+\ve)^2}$. Consider
$\gn(q(\ve),\ve-\o n)$. We have the simple estimates
\[ \lb{N1x}
\gn(q(\ve),\ve-\o n)\le \gn(q(\ve),-\o n)\le \gn(q_o,-\o n),
\]
since $q(\ve)\le q_o$ for $\ve <0$ small enough. Here  $\gn(q,s)$ is
the number of eigenvalues of $h_o-q$ on the set $(-\iy,s)$ counted
according to their multiplicity.
 \BBox


We estimate the number of eigenvalues of $\wt h$ on  the main gap
$\g=(-s_\g,0)$.

\no {\bf Proof of Theorem \ref{T1}.} i) Let $\gc={16\o\/s_e^2}$ and
$\d={1\/\gc}$. Define intervals $\g_-=[\m_0,0)$, where
$\m_0=-{s_\g\/2}$ and $I_j=[\m_{j-1},\m_j),\ \ j\in \N_{m},$ where
$$
 \m_j=\m_0+j\d, \ \ j\in \N_{m-1},\qq \m_m=0, \qq m-1=[{|\m_0|/\d}]\le
{8\o \/s_\g},
$$
where $[z]$ is the integer part of the number $z\ge 0$.

Due to $\g_-=\cup_1^m I_j$  we have $ \wt \gn(\g_-)=\sum \wt
\gn_j(I_j) $. Theorem \ref{Tei1} and  \er{N1}  yields
$$
\wt \gn(\g_-)=\sum _{j=1}^m\wt \gn(I_j),
 \qqq
\wt \gn(I_j)\le \sum _{n\in\Z}\gn_-(q_-,\m_j-\o |n|),\qqq
q_-=\gc u^2+v_-.
$$
We show \er{11}. We deduce that if $\o n>\|q_o\|_\iy$, then
$\gn_-(q_-,-\o n)=0$. Thus from Theorem \ref{Tks1}  we have
$\gn_-(q_-,-\o n)\le C_{ks} \|q_-\|_p^p$ for all $0\le n \o\le
\|q_-\|_\iy$, which yields
$$
\wt \gn(I_j)\le \sum _{n\in \Z}\gn_-(q_j,-\o |n|)\le C_{ks}
\|q_-\|_p^p \rt(1+{\|q_-\|_\iy\/\o}\rt),
$$
and then
$$
\wt \gn(\g_-)=\sum_1^m \wt \gn_j(I_j)\le m C_{ks} \|q_-\|_p^p
\rt(1+{\|q_-\|_\iy\/\o}\rt)<{9\o\/s_\g}C_{ks} \|q_-\|_p^p
\rt(1+{\|q_-\|_\iy\/\o}\rt),
$$
since $m\le 1+{8\o\/s_\g}<{9\o\/s_\g}$.

ii) We transform the case $\g_+$ to the case $\g_-$. We change the
variable $t'=-t$ and let $\pa'=-i{\pa\/\pa t'}$. Then the operator
$\wt h$ is unitary equivalent to the following operator
\[
\lb{U} \mU^* \wt h \mU= -\pa'+h_o+V_x(-t')=-\wt h', \qq \where \qq
\wt h'=(\pa'-h_o-V_x(-t')).
\]
for some unitary  operator $\mU$.  The spectrum of $\wt
h_o'=\pa'-h_o$ has the form
$$
\s(\wt h_o')=\bigcup_{n\in \Z} \s_n',\qqq   \s_n'= \s(-h_o)+\o n,\qq
\s_0'=\s(-h_o)\ss [-s_+,0].
$$
Since the spectrum of $\wt h'$ is $\o$ periodic, then we consider
the spectral interval $[-\o,0]$. We discuss the discrete spectrum of
$\wt h'$ in its main gap $\g'=(-\o,-s_+)$. The case
$\g_-'=[-{\o+s_+\/2},-s_+ )$ has been considered in i). Then we have
\er{12}.

iii) Collecting \er{11}  and \er{12} we obtain \er{12xx}. \BBox

We discuss estimates for $\o$ large  enough.

{\bf Proof of Theorem \ref{T2}.} i) Note that \er{14}, i.e., $\o\ge
s_++\d_-+\sqrt{\d_-(2s_++\d_-)}$, is equivalent to the following
inequality $\d_-c\le 1$, since $s_\g=\o-s_+$. Then  $\d_-$ and $c$
satisfy conditions from Theorem \ref{Tei1} ii) and we obtain
\[
\lb{N2xc} \wt \gn(I_-)\le \sum _{n\in \Z}\gn_-(q_-,-\o |n|).
\]
In order to show \er{15} we need to show that $q_-<\o$, which yields
$\gn(q_-,-\o n)=0$ for any $n\ge 1$. Using $c>{4\o\/s_\g^2}$ and
$\d_-c\le 1$ and $\d_+c\le \sqrt3$ we have two cases:

\no 1) Let $\d_+\le \d_-$. Then due to $c>{4\o\/s_\g^2}$ and $\d_-c\le
1$ we obtain
$$
q_-=c u^2+v_-\le c\d_-^2+\d_-=\d_-(c\d_-+1)\le 2\d_-\le
{2\/c}<{s_\g^2\/2\o}<{\o\/2}.
$$
2) Let $\d_+> \d_-$. Then due to $c>{4\o\/s_\g^2}$ and $\d_-c\le 1$
we obtain
$$
q_-=c u^2+v_-\le c\d_+^2+\d_-= {(c\d_+)^2\/c}+{c\d_-\/c}<{1\/c}(3+1)
\le {s_\g^2\/\o}<\o.
$$
This and \er{N2xc} yield $\wt \gn(I_-)\le \gn(q_-,0)$. Finally
Theorem \ref{Tks1} gives \er{15}.

ii) Using the transformation \er{U} and results of i) we obtain ii)

 iii) Summing  \er{14} and  \er{15}   we have
 \er{12xx}.
\BBox

\section {Eigenvalues on a.c. spectrum, proof of Theorem \ref{T3}}
\setcounter{equation}{0}

\subsection{Preliminaries} We consider the case $h_o=\D$ on $\Z^d, d\ge
1$.  Recall the  Rellich type theorem of Isozaki and Marioka
\cite{IM14}.

\begin{theorem}
\lb{TIM}

 Let $\l\in  (0, 2d)$ and $\gr_o>0$.  Suppose that a sequence
$(u_x)_{x\in \Z^d}$ satisfies $ (\D-\l)u=0$ in the domain $\{x\in
\Z^d: |x|>\gr_o\}$ and $\sum _{\gr_o<|x|<r}|u_x|^2=o(r)$ as $r\to
\iy$. Then there exists $\gr_1>\gr_o$ such that $f=0$ in $\{x\in
\Z^d: |x|>\gr_1\}$.
\end{theorem}

In this theorem if $\gr_o$ is fix, then there is no information
about $\gr_1$.

Below we need the  Rellich type theorem of Vesalainen \cite{V14}. It
is generalization of Theorem \ref{TIM} will concern vanishing in a
cone-like domain given by
\[
\lb{Km} \cK(\gm) =\{x\in \Z^d:
|x_1-\gm_1|+....+|x_{d-1}-\gm_{d-1}|\le (x_d-\gm_{d})\}.
\]

\begin{theorem} \lb{TV}


Let $u: Z^d\to \C $  be such that $\sum _{|x|<r}|u_x|^2=o(r)$  as \
$ r\to \iy$. And let $e^{t|x|}f_x\in \ell^2(\Z^d)$  for all $t>0$
and $ f\big|_{\cK(0)}=0$ for the cone-like domain given by \er{Km}.
Finally, let $\l\in  (0, 2d)$ and assume that $ (\D-\l)u=f$ in
$\Z^d$. Then $u\big|_{\cK(0)}=0$.

\end{theorem}

\no {\bf Proof of Theorem \ref{TKo}.} We consider the case $d=3$,
the proof of other cases is similar. Recall that the cone-like
domain $\cK(\gm), \gm\in \Z^3$ is defined by \er{Km}. Here $\gm\in
\Z^3$ is the top of the domain $\cK(\gm)$.

i) The half space $\dY_3(z), z\in \N$ is the union of cone-like
domains:
$$
\dY_3(z)=\bigcup_{\gn\in \Z^2} \cK(\gn,z+1).
$$
If $f |_{\dY_3(z)}=0$, then we obtain $f |_{\cK(\gn,z+1)}=0$ for
each $\gn\in \Z^2$. Then the  Vesalainen Theorem \ref{TV} gives $u
|_{\cK(\gn,z+1)}=0$ and then $u |_{\dY_3(z)}=0$.

ii) Let $\dF=\Z^3\sm \dQ(s)$ for some cuboid $\dQ(s), s=(s_j)\in
\N^3$. The domain $\dF$ is the union of the half space given by
$$
\begin{aligned}
\dF=\bigcup_{j=1,2,3} \dY_j^+ (s_j+1)\cup \dY_j^- (-1),
\end{aligned}
$$
where $\dY_j^\pm (z)=\{x\in \Z^3:  \pm x_j\ge \pm z\}$ is the half
space for $j=1,2,3$ and  some $z\in \N$.

If $f |_{\dF}=0$, then we obtain $f |_{\dY_j^+ (s_j+1)}=0$ and $f
|_{\dY_j^- (-1)}=0$ for all $j=1,2,3$.  The point i) gives $u
|_{\dY_j^+ (s_j+1)}=0$ and $u |_{\dY_j^- (-1)}=0$ for all $j=1,2,3$,
which yields $u |_{\dF}=0$. \BBox


Let $\p_j(x,t), j=1,..,N$ be eigenfunctions of $\wt h=\wt h_o+V$
such that
\[
\lb{ei1e} \wt h \p_j=\l_j \p_j,\qq \l_j\in (-\o,0).
\]
 Let $\gM(\m) $ be linear span of these eigenvectors $\p_j$.

In particular, we rewrite  the equation \er{ei1e} via the Fourier
transform by
\[
\lb{ei1e1}
\begin{aligned}
(h_o-(\l_j-\o n))\p_{j,n}+(V \p_j)_n=0, \qqq n\in \Z.
\end{aligned}
\]

\begin{lemma}
\lb{Tt3} Let $V$ satisfy \er{Vc} and $\c=\c_{\dQ(\gm)}$, where
$\dQ(\gm), \ \gm\in \N^d$ is the cuboid from \er{Vc}. Let
$\p=(\p_{n})_{n\in \Z}$ be an eigenfunction of $\wt h$, and $\wt h
\p=\l \p$. Let . Then it satisfies
\[
\lb{ee1} \c \p_{n}=\p_{n},\qq \forall \ n\in \N_{\gr}^-,\qq if \
\l\in \s_e=(-s_e,0),
\]
\[
\lb{ee2} \c \p_{n}=\p_{n},\qq \forall \ n\in \N_{\gr+1}^-, \qq if \
\l\in \s_b=(-\o,-s_e).
\]

\end{lemma}
\no{\bf Proof.}
 We rewrite the equation $\wt h \p=\l \p$  via the Fourier transformation:
$$
 0=(h_o-(\l_j-\o n )\p_{n}(x)+(V\p)_n(x)=(h_o-(\l_j-\o n )\p_{n}(x),
 \qqq \forall \ \ x\notin \dQ(\gm),
$$
 since $V$ satisfies \er{Vc}. Then the Rellich
type  Theorem \ref{TKo} yields \er{ee1} and \er{ee2}.
 \BBox


\subsection{Proof of Theorem \ref{T3}}
We discuss eigenvalues imbedded on the  interval $\s_b=(-\o,-s_e)$
of the continuous spectrum. We assume that each $V(t), t\in \T_\t$
is compactly supported. Define a potential $\vr_n^\gr$ by
\[
\lb{rn}
\vr_n^\gr=\ca \o |n|\c, \ & \ n\in \N_\gr^-\\
            0   , \ & \ n\notin \N_\gr^-\ac, \qq n\in \Z,\qq \where
            \  \gr=[s_+/\o]\in \Z.
\]
 where  $\c=\c_{\dQ(\gm)}$ and $\dQ(\gm), \ \gm\in \N^d$ is the cuboid from
\er{Vc}.

\begin{theorem}
\lb{Teac1x}

Let $\wt h=\wt h_o+V$, where $V$ satisfies \er{Vc}, \er{V1} and
$q_{-}=\gc_b u^2+v_-$. Let
\[
I_\d=(\m-\d,\m)\ss \s_b,\qqq \gc_b={4\o(\gr+2)\/s_e^2},\qq \d\le
{1\/\gc_b},\qq \gr=[s_+/\o].
\]
i) Then  the operator $\wt h$ has $\wt \gn(I_\d)$ eigenvalues on
$I_\d$ such that
\[
\lb{N1cq} \wt \gn(I_\d)\le \sum _{n\in
\Z}\gn_-(q_{-}+\vr_n^{\gr+1},\m-\o |n|),\qqq (q_--\o |n|+\m)_+.
\]
\no ii)  Then the operator $\wt h$ has $\wt \gn(\s_b)$ eigenvalues
on $\s_b$ such that
\[
\lb{132}
\begin{aligned}
 \wt \gn(\s_b)\le  C_b \rt(\gr+2+
2{\|q_-\|_\iy\/\o}\rt)C_{ks}^-\|q_-\|_p^p,
\end{aligned}
\]
where $C_b=1+  {4\o(\o-s_e)(\gr+2)\/s_e^2}$.

\end{theorem}

\no {\bf Proof.} i) Let $\p=\sum g_j\p_j$, where $\p, \p_{j}\in
\gM(I_\d)$. We need to estimate $\|\pa_-^{1\/2}\p\|^2$. Let
$P=P(\gr+1)$.  Using Lemma \ref{Tt3} we obtain
$$
(\pa_- P\p,\p)= (\pa_- P\c\p,\c\p)= \|\pa_-^{1\/2}P\c\p\|^2\le
(\gr+1)\o \|\c\p\|^2.
$$
Thus since there are  identities
\[
\|\pa_-^{1\/2}\p\|^2=\|\pa_-^{1\/2}P\p\|^2+\|\pa_-^{1\/2}P^\bot\p\|^2=
\|\pa_-^{1\/2}P\c \p\|^2+\|\pa_-^{1\/2}P^\bot\p\|^2,
\]
we need to estimate $\|\pa_-^{1\/2}P^\bot\p\|^2$. Using the identity
$$
\begin{aligned}
P^\bot\p =P^\bot R_o(\m)(\wt h_o-\m)\p=P^\bot R_o(\m)\big((\wt
h-\m)-V\big)\p,
\end{aligned}
$$
and the estimates $\|\pa_-^{1\/2}P^\bot R_o(\m)\|^2\le {c_b\/4}$ from
\er{r2} we have
\[
\lb{qwe1z}
\begin{aligned}
2\|\pa_-^{1\/2}P^\bot\p\|^2= 2\|\pa_-^{1\/2} P^\bot R_o(\m)\big((\wt
h-\m)-V\big)\p\|^2\le (\gc_b/2) \|\big((\wt h-\m)-V\big)\p\|^2
\\
\le \gc_b\big(\|(\wt h-\m)\p\|^2+\|V\p\|^2 \big).
\end{aligned}
\]
Let $T_o=\pa_++\pa_-+h_o$.  We have the identities
\[
\lb{qw1}
\begin{aligned}
((T_o-\m)\p,\p)=2(\pa_-\p,\p)-(V\p,\p)+((\wt h-\m)\p,\p)
\\
=2\|\pa_-^{1\/2}P\c\p\|^2+2(\pa_-P^\bot\p,\p) -(V\p,\p)+((\wt
h-\m)\p,\p).
\end{aligned}
\]
Then due to \er{qwe1z}  we obtain
$$
\begin{aligned}
((T_o-\m)\p,\p)\le
2\|\pa_-^{1\/2}P\c\p\|^2+\gc_b\Big(\|V\p\|^2+\|(\wt
h-\m)\p\|^2\Big)-(V\p,\p)+((\wt h-\m)\p,\p)
\\
\le ((\pa_-P\c+\gc_b u^2+v_-)\p,\p)= (q_-\p,\p),
\end{aligned}
$$
where we have used $\gc_b\|(\wt h-\m)\p\|^2+((\wt h-\m)\p,\p)\le 0$,
since $c_b(\l_j-\m)^2+(\l_j-\m)\le 0$. Then
$$
\begin{aligned}
\sum ((\D+\o |n|-\m)\p_n,\p_n)\le (q_-\p,\p)+\sum_{-r-1}^{-1} (\o \c
|n|\p_n,\p_n)=\sum_{\Z}((q_-+\vr_n^{\gr+1})\p_n,\p_n).
\end{aligned}
$$
ii) Let $\m_0\in
\s_b=(-\o, -s_e)$ and $m-1={\rm int} \{{|\m_m-\m_0|\/\d}\}$ is the
integer part of the number ${|\m_m-\m_0|\/\d}<{\o-s_e\/\d}$. Define
intervals $\s_b(\m)=[\m_0,\m_m)$ and
$$
 I_j=[\m_{j-1},\m_j),\ \ j\in \N_{m}, \qq \where\ \
 \m_j=\m_0+j\d, \ \ j\in \N_{m-1},\qq \m_m\in (\m_{m-1},-s_e).
$$
  The estimate \er{N1cq}  yields
$$
\wt \gn(I_j)\le \sum_{n\in \Z}\gn_-(q_{-}+\vr_n^{\gr+1},\m_j-\o |n|)
\le \cA, \qq \cA:=\sum_{n\in \Z}\gn_-(q_{-}+\vr_n^{\gr+1},-s_e-\o
|n|),
$$
since we take any numbers $\m_0\in \s_b=(-\o, -s_e)$ and $\m_m\in
(\m_{m-1},-s_e)$. Due to $\s_b(\m)=\cup_1^m I_j$  we deduce that
$$
\wt \gn(\s_b(\m))=\sum_1^m \wt \gn_j(I_j)\le m \cA<C_b \cA,
$$
where $C_b=1+{\o-s_e\/\d}=1+  {4\o(\o-s_e)(\gr+2)\/s_e^2}$. We have
$\cA=\cA_1+\cA_2$, where
$$
\cA_1=\sum_{n=-\gr-1}^{0}\gn_-(q_{-}+\c \o |n|,-s_e-\o |n|),\qq
\cA_2=\sum_{n\in G}\gn_-(q_{-},-s_e-\o |n|),
$$
and $G=\Z\sm \{-\gr-1,...,0\}$. Consider $\cA_1$. From \er{n-} we
obtain
$$
\gn_-(q_{-}+\c \o |n|,-s_e-\o |n|)\le \gn_-(q_{-},-s_e)\le
C_{ks}^-\|(q_--s_e)_+)\|_p^p,
$$
which yields
$$
\cA_1=\sum_{n=-\gr-1}^{0}\gn_-(q_{-}+\c \o |n|,-s_e-\o |n|)\le
(\gr+2)C_{ks}^-\|(q_--s_e)_+)\|_p^p.
$$
Consider $\cA_2$. We deduce that if $\o |n|>\|q_-\|_\iy$, then
$\gn_-(q_-,-s_e-\o |n|)=0$. Thus from Theorem \ref{Tks1}  we have
$\gn_-(q_-,-\o |n|)\le C_{ks}^- \|q_-\|_p^p$ for all $1\le |n| \o\le
\|q_-\|_\iy$, which yields
$$
\cA_2=\sum_{n\in G}\gn_-(q_{-},-s_e-\o |n|)\le
2{\|q_-\|_\iy\/\o}C_{ks}^- \|q_-\|_p^p.
$$
Collecting all estimates we obtain \er{132}. \BBox


We discuss eigenvalues on the interval $\s_e=(-s_e,0)$.

\begin{theorem}
\lb{Teac2}

Let $\wt h=\wt h_o+V$, where $V$ satisfies \er{Vc} and
$\gc_e={16\o(\gr+1)\/s_e^2}$.

i) Let $\m\in \s_e^-=[-{s_e\/2},0)$. Let $V$ satisfy \er{V1} and
$q_-=\gc_eu^2+v_-$.
 Then  the operator $\wt h$ has $\wt \gn(\s_e^-)$ eigenvalues on $\s_e^-$ such that
\[
\lb{N1qq} \wt \gn(\s_e^-)\le C_e \rt(\gr+1+
2{\|q_-\|_\iy\/\o}\rt)C_{ks}^-\|q_-\|_p^p,
\]
where $C_e=1+{8\o(\gr+1)\/s_e}$.

ii) Let $V$ satisfies \er{V2}, and a potential $q_+=c_e u^2+v_+$.
Let $\s_e^+=(-s_e,-{s_e\/2}]$.  Then
\[
\lb{N2c} \wt \gn(\s_e^+)\le C_e \rt(\gr+1+
2{\|q_+\|_\iy\/\o}\rt)C_{ks}^+\|q_+\|_p^p.
\]

\end{theorem}

\no {\bf Proof.} i) Let $I_\d=[\m-\d,\m)\ss \s_e^-$, where $\d\le
{1\/C}$. Let $\p=\sum c_j\p_j$, where $\p, \p_{j}\in \gM(I_\d)$. We
need to estimate $\|\pa_-^{1\/2}\p\|^2$. Let $P=P(\gr)$. We have
identities
\[
\|\pa_-^{1\/2}\p\|^2=\|\pa_-^{1\/2}P\p\|^2+\|\pa_-^{1\/2}P^\bot\p\|^2=
\|\pa_-^{1\/2}P\c \p\|^2+\|\pa_-^{1\/2}P^\bot\p\|^2.
\]
 Using Lemma \ref{Tt3} we obtain
$$
(\pa_- P\p,\p)= (\pa_- P\c\p,\c\p)=
\|\pa_-^{1\/2}P\c\p\|^2=\sum_{-\gr}^{-1} (\vr_n^{\gr}\p_n,\p_n).
$$
We need to estimate $\|\pa_-^{1\/2}P^\bot\p\|^2$. Using the estimate
$\|\pa_-^{1\/2}P^\bot R_o(\m)\|^2\le c_e/4$ from \er{r1} and the
identity
\[
\lb{12e}
\begin{aligned}
P^\bot\p =R_o(\m)P^\bot(\wt h_o-\m)\p=R_o(\m)P^\bot\big((\wt
h-\m)-V\big)\p,
\end{aligned}
\]
we have
\[
\lb{qwe1}
\begin{aligned}
2\|\pa_-^{1\/2}P^\bot\p\|^2= 2\|\pa_-^{1\/2} P^\bot R_o(\m)\big((\wt
h-\m)-V\big)\p\|^2\le (c_e/2) \|\big((\wt h-\m)-V\big)\p\|^2
\\
\le c_e\big(\|(\wt h-\m)\p\|^2+\|V\p\|^2 \big).
\end{aligned}
\]
Thus using the identities
$$
\begin{aligned}
((T_o-\m)\p,\p)=2(\pa_-\p,\p)-(V\p,\p)+((\wt h-\m)\p,\p)
\\
=2\|\pa_-^{1\/2}P\c\p\|^2+2(\pa_-P^\bot\p,\p) -(V\p,\p)+((\wt
h-\m)\p,\p)
\end{aligned}
$$
and the estimates \er{qwe1}, \er{qwe1}  we obtain
$$
\begin{aligned}
((T_o-\m)\p,\p)\le 2\|\pa_-^{1\/2}P\c\p\|^2+c_e\Big(\|V\p\|^2+\|(\wt
h-\m)\p\|^2\Big)-(V\p,\p)+((\wt h-\m)\p,\p)
\\
\le (q_-\p,\p)+\sum_{-r-1}^{-1} (\o \c
|n|\p_n,\p_n)=\sum_{\Z}(q_-+\vr_n^{\gr+1}\p_n,\p_n),
\end{aligned}
$$
where we have used $c_e\|(\wt h-\m)\p\|^2+((\wt h-\m)\p,\p)\le 0$,
since $c_e(\l_j-\m)^2+(\l_j-\m)\le 0$. This yields
 that  the operator $\wt h$ has $\wt \gn(I_\d)$ eigenvalues on $I_\d$ such that
\[
\lb{N1c} \wt \gn(I_\d)\le \sum _{n\in \Z}\gn_-(q_-+\vr_n^{\gr},\m-\o
n),\qqq (q_--\o n+\m)_+.
\]
 We show \er{N1qq}.  Let $\m_0=-{s_e\/2}, \d={1\/c_e}$ and $m-1={\rm int}
\{{|\m_m-\m_0|\/\d}\}$ is the integer part of the number
${|\m_m-\m_0|\/\d}<{s_e\/2\d}$. Define intervals
$\s_e^-(\m)=[\m_0,\m_m)$ and
$$
 I_j=[\m_{j-1},\m_j),\ \ j\in \N_{m}, \qq \where\ \
 \m_j=\m_0+j\d, \ \ j\in \N_{m-1},\qq \m_m\in (\m_{m-1},0).
$$
  The estimate \er{N1c}  yields
$$
\wt \gn(I_j)\le \sum_{n\in \Z}\gn_-(q_{-}+\vr_n^{\gr},\m_j-\o |n|)
\le \cB, \qq \cB:=\sum_{n\in \Z}\gn_-(q_{-}+\vr_n^{\gr},-\o |n|),
$$
where $q_-=Cu^2+v_-$, since we take any number $\m_m\in
(\m_{m-1},0)$. Due to $\s_e^-(\m)=\cup_1^m I_j$  we have
$$
\wt \gn(\s_e^-(\m))=\sum_1^m \wt \gn_j(I_j)\le m \cA<C_e \cB,
$$
where $C_e=1+{s_e\/2\d}=1+{8\o(\gr+1)\/s_e}$. We have
$\cB=\cB_1+\cB_2$, where
$$
\cB_1=\sum_{n=-\gr}^{0}\gn_-(q_{-}+\c \o |n|,-\o |n|),\qq
\cB_2=\sum_{n\in G}\gn_-(q_{-},-\o |n|),
$$
and $G=\Z\sm \{-\gr,...,0\}$. Consider $\cB_1$. From \er{n-} we
have
$$
\gn_-(q_{-}+\c \o |n|,-\o |n|)\le \gn_-(q_{-},0)\le
C_{ks}^-\|q_-)\|_p^p
$$
which yields
$$
\cB_1=\sum_{n=-\gr}^{0}\gn_-(q_{-}+\c \o |n|,-\o |n|)\le
(\gr+1)C_{ks}^-\|q_-\|_p^p.
$$
Consider $\cB_2$. We deduce that if $\o |n|>\|q_-\|_\iy$, then
$\gn_-(q_-,-\o |n|)=0$. Thus from Theorem \ref{Tks1}  we have
$\gn_-(q_-,-\o |n|)\le C_{ks}^- \|q_-\|_p^p$ for all $1\le |n| \o\le
\|q_-\|_\iy$, which yields
$$
\cB_2=\sum_{n\in G}\gn_-(q_{-},-\o |n|)\le
2{\|q_-\|_\iy\/\o}C_{ks}^- \|q_-\|_p^p.
$$
Collecting all estimates we obtain \er{132}.

ii) We transform the case $\s_e^+$ to the case $\s_e^-$. We change
the variable $t'=-t$ and let $\pa'=-i{\pa\/\pa t'}$. Then the
operator $\wt h$ is unitary equivalent to the following operator
\[
\lb{Ux} \mU^* \wt h \mU= -\pa'+h_o+V_x(-t')=-\wt h', \qq \where \qq
\wt h'=(\pa'-h_o-V_x(-t')),
\]
for some unitary  operator $\mU$.  The spectrum of $\wt
h_o'=\pa'-h_o$ has the form
$$
\s(\wt h_o')=\bigcup_{n\in \Z} \s_n',\qqq   \s_n'= \s(-h_o)+\o n,\qq
\s_0'=\s(-h_o)\ss [-s_+,0].
$$
Since the spectrum of $\wt h'$ is $\o$ periodic, then we consider
the spectral interval $[-\o,0]$. We discuss the discrete spectrum of
$\wt h'$ in its main gap $\g'=(-\o,-s_+)$. The case
$\gI_-'=[-(\o+s_+)/2,-s_+ )$ has been considered in i). Then we have
\er{N2c}. \BBox


\no {\bf Proof of Theorem \ref{T3}.} The inequality \er{19} have been
proved in Theorem \ref{Teac1x}. Summing two estimates of $\wt
\gn(\s_e^\pm)$ from
 Theorem \ref{Teac2} we obtain  \er{Nne}.
 \BBox

\section {Eigenvalues  for the 1dim case }
\setcounter{equation}{0}

\subsection {Stationary case } We  discuss
one-dimensional stationary case. We consider Schr\"odinger operators
$ (hy)_n=(h_oy)_n+\gq_ny_n$ on $\N$, where  $h_o=\D$ is the
Laplacian on $\N$ given by
\[
\lb{J1} (\D y)_n= y_{n-1}+y_{n+1}+2y_n,\qqq  n\in\N,
\]
and formally $y_0=0$. Assume that the potential $\gq=(\gq_n)_1^\iy$
is real and belongs to the space $\ell_1^1(\N)$ of sequences
$y=(y_n)_1^\iy$ equipped with the norm $\|y\|_{\bu}=\sum _{n\ge 1}
n|y_n|$. It is well known that the operator $h$ has purely
absolutely continuous spectrum $[0, 4]$ plus a finite number of
simple eigenvalues on the set $\R\sm[0,4]$. We recall results about
the discrete analogue of Bargmann's bound from \cite{G82},
\cite{G88}, \cite{HS02}.

\begin{theorem}
\lb{THS}

Let $h=h_o+\gq$ on $\N$, where $\gq\in \ell_1^1(\N)$. Let $\gn(\gq)$
be the number of eigenvalues of $h$ outside $[-2, 2]$. Then
\[
\lb{J3} \gn(\gq)\le \|\gq\|_{\bu}.
\]
Let $\gn_\pm(\gq)$ be the number of eigenvalues $\l$ with $\pm \l >
2$. Then the number $\gn_\pm(\gq)$ of eigenvalues of $h$ satisfies
\[
\lb{J4} \gn_\pm(\gq)\le \|(\gq)_\pm\|_{\bu},
\]
where $(z)_+=\max \{0,z\}, z\in \R$.

\end{theorem}

\subsection {Time periodic potentials } Consider the Schr\"odinger
equation on $\N$:
$$
{\tfrac{ d }{d t}}y(t)= -i h(t) y(t), \qq \where \ \ h(t)=h_o+ V(t),
$$
and $h(t)$  is the Hamiltonian, $\t$-periodic in time $t$ and the
operator $h_o=\D$ is defined by \er{J1}.  We assume that the
potential $V_x(t), (t,x)\in \R\ts \N$ is real
 $\t$-periodic in time   and satisfies
\[
\lb{J5} \sup_{t\in \T_\t}|V(t)|\le u,\qqq   \where \qq u^2\in
\ell_1^1(\N).
\]
Introduce the operators $\wt h_o=\pa+h_o$ and   $\wt h=\wt h_o+V$ on
$L^2(\T_\t,\ell^2(\N))$. Recall that the spectrum of $\wt h_o$ is
abs. continuous and has the form
$
 \s(\wt h_o)=\cup_{n\in \Z} \s_n$, where the band $\s_n=
[0,4]+\o n$. Note that if $\o >4$, then the spectrum of $\wt h_o$ is
a union of bands separated by gaps.

Firstly, we consider the case $\o>4$, when there is the main gap
$\g\ss (-\o,0]$ defined by
\[
\lb{hd}  \g=(-s_\g,0)\ne  \es, \qq \where \qq s_\g=\o-4>0.
\]
We estimate the number $\wt \gn(\g)$of eigenvalues of $\wt h$ on the
main gap $\g$.

\begin{corollary}
\lb{TJ1}

 Let $\wt h=\wt h_o+V$ on
$L^2(\T_\t,\ell^2(\N))$ and let      $s_\g=\o-4>0$.

\no i)  Let $\g^+=(-s_\g,-{s_\g\/2}]$ and $\g_-=[-{s_\g\/2},0)$.
Assume that potentials $V$ and $q_\pm=\gc u^2+v_\pm$, where
$\gc={16\o\/s_\g^2}$ satisfy
\[
\lb{J6}   \sup_{t\in \T_\t} V(t)_\pm\le v_\pm,\qqq \where \qq
v_\pm\in \ell_1^1(\N).
\]
 Then the operator $\wt h$ has $\wt \gn(\g_\pm)$ eigenvalues in $\g_\pm$
 such that
\[
\lb{J7}
\begin{aligned}
 \wt \gn(\g_\pm)\le
 {9\o\/s_\g}  \rt(1+2{\|q_\pm\|_\iy\/\o}\rt)\|q_\pm\|_\bu.
\end{aligned}
\]
ii) Let $V$ satisfy \er{J6}. Then the number of eigenvalues of the
operator $\wt J$ in the main gap $\g=(-s_\g,0)$ satisfies
\[
\lb{J8}
\begin{aligned}
 \wt \gn(\g)\le \sum _{\pm}{9\o\/s_\g}  \rt(1+2{\|q_\pm\|_\iy\/\o}\rt)\|q_\pm\|_\bu.
\end{aligned}
\]

\end{corollary}

\no {\bf Proof.} We omit  the proof, since it  is similar to the proof
of Theorem \ref{T1}.
 \BBox

We estimate the number of eigenvalues of $\wt h$ on the main gap
$\g$, when $\o$ and $s_\g=\o-4>0$ are large  enough.

 \begin{corollary}
\lb{TJ2}

Consider  $\wt h=\wt h_o+V$ on $L^2(\T_\t,\ell^2(\N))$. Let $V$
satisfy \er{J6} and let a potential $q_\pm=c u^2+v_\pm$, where
$$
 c>{4\o\/s_\g^2},\qq \d c\le \sqrt3,\qq  \o\ge
4+\d+\sqrt{\d(8+\d)},\qq \d=\sup_{(t,x)\in \T_\t\ts \N}
|V_x(t)|<{s_\g\/2}.
$$
 Then the
operator $\wt h$ has $\wt \gn(\g)$ eigenvalues on the main gap
$\g=(-s_\g,0)$ such that
\[
\lb{18c}
\begin{aligned}
 \wt \gn(\g)\le \|q_{_+}\|_\bu+\|q_{_-}\|_\bu.
\end{aligned}
\]

\end{corollary}

\no {\bf Proof.} We omit  the proof, since it  is similar to the proof
of Theorem \ref{T2}.
 \BBox

We consider the imbedding eigenvalues of $\wt h$.  The main spectral
interval has the form $(-\o,0)=\s_e\cup \s_b\cup \{s_e\}$, where
\[
\lb{J9} \gr=[4/\o]\ge 0,\qqq  s_e=\gr \o-4,\qq \s_b=(-\o, -s_e),\qq
\s_e=(-s_e,0).
\]
There are two cases:

\no if $\o>4$, i.e., $\gr=0$, then there are a gap $\g=(-s_b,0)$ and
a band $\s_b=(-\o, -s_b)$ on $(-\o,0)$,

\no if  $\o\le 4$,  i.e., $\gr\ge 1$, then $\s(\wt h_o)=\R$ and
there are no gaps in the spectrum of $\wt h_o$.

We discuss the number of eigenvalues imbedding on the main spectral
interval $(-\o,0)$.

\begin{corollary}
\lb{TJ3}

Let $\wt h=\wt h_o+V$ on $L^2(\T_\t,\ell^2(\Z^d))$. Assume that $V$
satisfy \er{J6}  and $ \supp V\ss \T_\t\ts [1,...,K]$ for some $K\in
\N$. Define  potentials
$$
q_{\a,\pm}=\gc_\a u^2+v_\pm,\qq  \a\in \{b,e\},\qq \gc_e=4\gc_b,\qq
\gc_b={4\o(\gr+2)\/s_e^2},\qq \gr=[4/\o]\ge 0.
$$

\no i) Let $\gr\ge  0$. Then the operator $\wt h$ has $\wt
\gn(\s_b)$ eigenvalues on $\s_b=(-\o, -s_e)$ such that
\[
\lb{19c} \wt \gn(\s_b)\le  C_b \rt(\gr+2+
2{\|q_{b,\pm}\|_\iy\/\o}\rt)\|q_{b,\pm}\|_\bu,
\]
where $C_b=1+  {4\o(\o-s_e)(\gr+2)\/s_e^2}$.

\no ii) Let  $\gr\ge 1$. Then the operator $\wt h$ has $\wt
\gn(\s_e)$ eigenvalues on $\s_e=(-s_e, 0), s_e=\gr \o-4$ such that
\[
\lb{Nnec} \wt \gn(\s_e)\le C^+\|q_{e,+}\|_\bu+C^-\|q_{e,-}\|_\bu,
\]
where $ C^\pm=C_e \rt(\gr+1+ 2{\|q_{e,\pm}\|_\iy\/\o}\rt)$ and
$C_e=1+{8\o(\gr+1)\/s_e}$.

\end{corollary}

\no {\bf Proof.} We omit  the proof, since it  is similar to the proof
of Theorem \ref{T3}.
 \BBox

\

\footnotesize\footnotesize
  \no {\bf Acknowledgments.}
  E. K. was supported by the RSF grant  No.
19-71-30002.


\end{document}